\documentclass[12pt]{amsart}
\usepackage{amsfonts,amsmath,amsthm,amscd,amssymb,latexsym,cite,verbatim,texdraw,floatflt,
caption2}
\RequirePackage{hyperref}

\usepackage[a4paper]{geometry}

\theoremstyle
{plain}

\begin{document}

\title{Sequential coarse structures of topological groups}

\author{Igor Protasov}

\maketitle
\vskip 5pt

{\bf Abstract.}
We endow a topological group $(G, \tau)$  with a coarse structure defined by the smallest group ideal
$S_{\tau} $ on $G$ containing all  converging sequences  with their limits and denote the obtained coarse group by
 $(G, S_{\tau})$. If $G$  is discrete then
 $(G, S_{\tau})$
  is a finitary coarse group studding in  {\it Geometric Group Theory}.
The main result: if a topological abelian group
$(G, \tau)$
 contains a non-trivial converging sequence then {\it asdim}  $(G, S_{\tau})= \infty $.
 We study metrizability, normality and functional boundedness of  sequential coarse groups and put some open questions.
\vskip 10pt

\vskip 10pt


\section{Introduction}

Let $X$  be a set. A family $\mathcal{E}$ of subsets of $X\times X$ is called a {\it coarse structure } if
\vskip 7pt

\begin{itemize}
\item{}   each $E\in \mathcal{E}$  contains the diagonal  $\bigtriangleup _{X}$,
$\bigtriangleup _{X}= \{(x,x): x\in X\}$;
\vskip 5pt

\item{}  if  $E$, $E^{\prime} \in \mathcal{E}$ then $E\circ E^{\prime}\in\mathcal{E}$ and
$E^{-1}\in \mathcal{E}$,   where    $E\circ E^{\prime}=\{(x,y): \exists z((x,z) \in  E,  \   \ (z, y)\in E^{\prime})\}$,   $E^{-1}=\{(y,x): (x,y)\in E\}$;
\vskip 5pt

\item{} if $E\in\mathcal{E}$ and $\bigtriangleup_{X}\subseteq E^{\prime}\subseteq E  $   then
$E^{\prime}\in \mathcal{E}$;
\vskip 5pt

\item{}  for any   $x,y\in X$, there exists $E\in \mathcal{E}$   such that $(x,y)\in E$.

\end{itemize}
\vskip 7pt

A subset $\mathcal{E}^{\prime} \subseteq \mathcal{E}$  is called a
{\it base} for $\mathcal{E}$  if, for every $E\in \mathcal{E}$, there exists
  $E^{\prime}\in \mathcal{E}^{\prime}$  such  that
  $E\subseteq E ^{\prime}$.
For $x\in X$,  $A\subseteq  X$  and
$E\in \mathcal{E}$, we denote
$E[x]= \{y\in X: (x,y) \in E\}$,
 $E [A] = \cup_{a\in A}   \   \   E[a]$
 and say that  $E[x]$
  and $E[A]$
   are {\it balls of radius $E$
   around} $x$  and $A$.

The pair $(X,\mathcal{E})$ is called a {\it coarse space} \cite{b14}  or a ballean \cite{b9}, \cite{b12}.

Each subset $Y\subseteq X$  defines the {\it subballean}
$(Y, \mathcal{E}_{Y})$, where $\mathcal{E}_{Y}$  is the restriction of $\mathcal{E}$  to $Y\times Y$.
A subset $Y$  is called {\it bounded} if $Y\subseteq E[x]$  for some $x\in X$  and
$E\in \mathcal{E}$.
\vskip 7pt

A family $\mathcal{F}$  of subsets  of $X$  is called $E$-{\it bounded ($E$-disjoint)} if, for each $A\in  \mathcal{F}$, there exists $x\in X$ such that $A\subseteq  E[x]$   $(E[A]\cap   E[B] = \emptyset$  for all distinct $A, B \in \mathcal{F})$.

By the definition [14, Chapter 9],  {\it asdim} $(X, \mathcal{E})\geq  n$  if, for each $E\in  \mathcal{E}$,   there exist $F\in \mathcal{E}$  and $F$-bounded covering $\mathcal{M}$ of $X$  which can be partitioned $\mathcal{M}= \mathcal{M} _{0}\cup  \ldots \cup \mathcal{M}_{n}$
 so that each family $M_{i}$ is $E$-disjoint.
If there exists the minimal $n$  with this property then
 {\it asdim }$(X, \mathcal{E})= n$, otherwise  {\it asdim } $(X, \mathcal{E}) = \infty $.

Given two coarse spaces  $(X, \mathcal{E})$,  $(X^{\prime}, \mathcal{E}^{\prime})$, a mapping
   $f: X\longrightarrow  X^{\prime}$   is called {\it macro-uniform}  if, for each
    $E \in \mathcal{E}$,  there  exists
    $E^{\prime} \in \mathcal{E}^{\prime}$
     such that
     $f(E[x]) \subseteq  E^{\prime}[f(x)]$ for each $x\in X$.
If $f$  is a bijection such that $f$  and $f ^{-1}$  are macro-uniform  then $f$  is called an {\it asymorphism}.

Now let $G$  be a group. A family $\mathcal{I}$  of  subsets of $G$  is
 called a {\it group ideal} \cite{b10}, \cite{b12}  if $G$  contains the family
  $[G] ^{< \omega}$  of all finite subsets of $G$  and
  $A, B \in \mathcal{I}$,
   $ \  \   C\subseteq A$
    imply
   $AB ^{-1}  \in \mathcal{I}$,  $C\in \mathcal{I}$.
Every group ideal $\mathcal{I}$  defines a coarse structure on $G $
  with the base $\{\{ (x, y): x\in Ay \} : A \in \mathcal{I}\}$.
We denote $G$  endowed with this coarse structure by  $(G, \mathcal{I})$.

If $G$  is discrete  then the coarse space
 $(G, [G] ^{< \omega })$  is the main subject of {\it Geometric Group Theory}, see \cite{b5}.
For coarse structures on $G$  defined by the ideal $[G] ^{< \kappa} $, where $\kappa$  is a cardinal, see \cite{b11}.

Every topological group $G$  can be endowed with a coarse structure defined by the ideal of all totally bounded subsets of  $G$.
These coarse structures were introduced and studied in \cite{b6}.
For asymptotic dimensions of locally compact abelian groups endowed with coarse structures defined by ideals of  precompact subsets see \cite{b7}.
We recall that a subset $Y$  of a topological space is {\it precompact} if the closure  of $Y$   is  compact.

For a topological group
$(G, \tau)$,
 we denote by  $\mathcal{C}_{\tau}$  the group ideal of precompact  subsets of $G$,
  and by
   $\mathcal{S}_{\tau}$
    the minimal group ideal containing  all  converging sequences with their limits.
    Clearly,
 $\mathcal{S}_{\tau}\subseteq \mathcal{C}_{\tau}$.
\vskip 10pt

\section{Asympotic dimension}

We  recall \cite{b13} that a sequence $(a_{n}) _{n \in\omega} $
 in an abelian group  $G$ is a {\it $T$-sequence} if there exists a Hausdorff group  topology on
 $G$  in which $(a_{n}) _{n \in\omega} $  converges to $0$.
For  a  $T$-sequences
  $(a_{n}) _{n \in\omega} $ on
   $G$, we denote by  $\tau _{(a_{n})}$  the strongest group topology on $G$
in which $(a_{n}) _{n \in\omega} $
converges to $0$.
We put $A=\{ 0,  a_{n}, -a_{n}: n \in \omega\}$  and denote by
$A_{n}$  the sum on $n$  copies of $A$.

By [13, Theorem 2.3.11], $(G, \tau _{(a_{n})})$   is complete.
Hence, a subset $S$  of $G$  is totally bounded in $(G, \tau _{(a_{n})})$ if  and only if $S$  is precompact.

We use the following three theorems proved by the author in \cite{b3}.

\vspace{7 mm}

{\bf Theorem 1}. {\it For any $T$-sequences    $(a_{n}) _{n \in\omega} $ on $G$,  the family
    $\{F+ A_{n}:  F\in [G]^{<\omega}\}$, $n\in\omega\}$
     is  a base for the ideal
           $\mathcal{C}_{\tau(a_{n})}$
            and
            $\mathcal{C}_{\tau(a_{n})}=\mathcal{S}_{\tau(a_{n})}$.
         If  $G $  is generated by the set
      $\{a_{n}: n\in\omega \}$ then $\{A_{n} : n\in\omega\}$
        is a base for
          $\mathcal{C}_{\tau(a_{n})}$.
\vspace{5 mm}

Proof.}
Apply Lemma 2.3.2 from \cite{b13}.  \hfill  $\Box$

 \vspace{8 mm}

Given an arbitrary subset  $S$ of $G$,  the Cayley  graph Cay $(G, S)$ is a graph with the set of vertices  $G$  and the set of edges  $\{(x, y): x-y\in S\cup (-S)\}$.  \vspace{7 mm}

{\bf Theorem 2}. {\it If a
$T$-sequences    $(a_{n}) _{n \in\omega} $
 generates $G$  then the coarse group
 $(G, \mathcal{C}_{\tau(a_{n})})$
 is asymorphic  to
Cay $(G,\{a_{n}: n\in\omega \})$.

\vspace{5 mm}

Proof.}
Apply and Theorem  1 and  Theorem   5.1.1 from \cite{b12}. \hfill  $\Box$

\vspace{6 mm}

We recall  that  the Hamming space $\mathbb{H}$ is the set
$[\omega ]^{<\omega} $   endowed with the metric
$h(F, H) = |F \bigtriangleup H|$.
To see that
{\it asdim}  $ \mathbb{H} = \infty $,
 it suffices to find  an asymorphic copy of $\mathbb{N}$  in $\mathbb{H}$  and observe that $\mathbb{H}$ is asymorphic to $\mathbb{H}\times \mathbb{H}$.

 \vspace{5 mm}

{\bf Example 1.}
Let $G$  be the direct sum of groups
$\{\langle a_{n}\rangle : n\in\omega\}$
 of order 2.
Clearly,
$(a_{n}) _{n\in\omega} $
 is a
 $T$-sequence
 on $G$.
By Theorem 1,
 the canonical bijection between
 $( G, \mathcal{C}_{\tau(a_{n})})$
  and the Hamming space $\mathbb{H}$  of all finite subsets of $\omega$  is an asymorphism.

A  $T$-sequences    $(a_{n}) _{n \in\omega} $
is  called {\it trivial} if  $a_{n}= 0$  for all but finitely many $n\in \omega$.

\vspace{7 mm}

{\bf Theorem 3}. {\it
For any non-trivial $T$-sequences    $(a_{n}) _{n \in\omega} $
  on $G$, the coarse group
   $( G, \mathcal{C}_{\tau(a_{n})})$
  contains a subspace asymorphic to the Hamming space $\mathbb{H}$ so
 {\it asdim} $ ( G, \mathcal{C}_{\tau(a_{n})}) = \infty$.

\vspace{7 mm}

Proof.}
Without loss of generality,  we suppose that
$\{a_{n}: n\in\omega\}$
 generates $G$  and $a_{n}\neq 0$  for each $n\in \omega$.

Given an arbitrary $T$-sequence    $(b_{n}) _{n \in\omega} $ in  $G$, we denote
   $$
   FS (b_{n}) _{n \in\omega} = \left\{\sum_{i\in F} b _{i}:  F\in [\omega]^{<\omega} \right\}
   $$
   and say that   $(b_{n}) _{n \in\omega} $
    is  FS-{\it strict} if, for any
    $H, F\in [\omega]^{<\omega}$,
    $$
    \sum_{i\in H}  b_{i}= \sum_{i\in F}  b_{i} \  \  \  \Longrightarrow \  \  \  F=H.
    $$

   We note that
$(b_{n}) _{n \in\omega}$
is  FS-strict  if, for each  $n\in\omega $,
\vspace{1mm}
$$
 b_{n+1}\notin \left\{\sum_{i\in F}  b_{i} -  \sum_{i\in H}  b_{i}: H, F\subseteq\{0, \ldots , n\} \right\}. \eqno(1)
$$
\vspace{1 mm}

We   assume that $(b_{n}) _{n \in\omega}$ is FS-strict and
\vspace{1 mm}

$$
 \mbox{if } \ b=\sum_{i\in F}  b_{i}, \    a\in A_{n}
 \mbox{ and } b+a=\sum_{i\in H}  b_{i}   \mbox{ then } |F \triangle H|\leq n. \eqno(2)
$$
\vspace{1mm}

Then the canonical bijection $f: \mathbb{H}\longrightarrow FS (b_{n}), \  \  f(H)= \sum_{i\in H}  b_{i}$ is an asymorphism.

To construct the desired sequence  $(b_{n})_{n\in\omega}$ we rewrite  $(2)$
 in the following equivalent form
\vspace{1 mm}

 $$
  \mbox{ if } \ i_{0} <  i_{1} <  \ldots < i_{n} <  \omega \  \mbox{ and } t _{i_{0}} ,    \ldots , t _{i_{n}} \in  \{1, -1\}
  \mbox{ then } t_{i_{0}} \  b _{i_{0}}+     \ldots +  t _{i_{n}} b _{i_{n}}\notin  A_{n}.   \eqno(3)
$$

\vspace{1mm}

We put $b_{0}= a_{0}$ and assume that  $b_{0}, \ldots , b_{n}$  have been chosen. We show how to choose  $b _{n+1}$ to satisfy $(1)$  and
\vspace{1 mm}

  $$
  \sum_{s=0}^kt _{i_{s}} b _{i_{s}}
 + t _{n+1} b _{n+1} \notin  A_{k+1}\mbox{ for }i_{0}< \ldots < i_{k}\leq n
 \mbox{ and } $$  $$t_i\in\{1, -1\}\mbox{ for } i \in \{i_0,\ldots,i_k, {n+1}\}.\eqno(4)
$$

  \vspace{2 mm}

We  assume that there exists a subsequence $(c_{m})_{m\in\omega}$
of $(a_{n})_{n\in\omega}$
such that $t _{i_{0}} b _{i_{0}} + \ldots +  t _{i_{k}} b _{i_{k}} +
 t  c _{m} \in  A_{k+1}$ for $t \in  \{1, -1\} $ and for each $m\in\omega$.
Every infinite subset of $A _{k+1}$  has a limit point in $A _{k}$.
Hence, $t _{i_{0}} b _{i_{0}} + \ldots +  t _{i_{k}} b _{i_{k}}  \in  A_{k+1}$
 contradicting the  choice of $b_{0}, \ldots , b_{n}$. Thus, $b _{n+1}$  can be taken from
$\{a_{m_{0}}, a_{m_{0}+1},  \ldots\}$ for some $m_{0} \in \omega$.  \hfill  $\Box$

\vspace{6 mm}

{\bf Theorem 4}.  {\it
Let  $G$ be a group and let $\mathcal{I}$  be a group ideal on $G$.
Assume that  there exists a family  $\mathcal{F}$ of group ideals on $G$  such that
$\mathcal{J} \subseteq \mathcal{I}$ for each  $\mathcal{J}\in \mathcal{F}$  and, for each $A\in \mathcal{I}$,  there exists
$\mathcal{J}\in \mathcal{F}$   such that $A\in \mathcal{J}$.
If {\it asdim} $(G, \mathcal{J}) \geq n $  for each  $\mathcal{J}\in \mathcal{F}$   then
{\it asdim} $ (G, \mathcal{I}) \geq n $.

\vskip 7pt

Proof.}
Given any $A\in \mathcal{I} $  and $n\in \omega$,  we choose
$\mathcal{J}\in \mathcal{F}$
 such that
 $A\in \mathcal{J}$
  and a uniformly bounded covering $\mathcal{M}$ of $(G, \mathcal{J})$  and a partition
   $\mathcal{M} _{0},\ldots , \mathcal{M}_{n}$ of $\mathcal{M}$  witnessing
{\it asdim}  $ (G, \mathcal{J}) \geq n $.
Then these
$\mathcal{M}, \mathcal{M} _{0},\ldots , \mathcal{M}_{n}$
 say that
 {\it asdim} $ (G, \mathcal{I}) \geq n $.  \hfill   $    \Box$

\vspace{6 mm}

{\bf Theorem 5}.  {\it
If a topological abelian group $(G, \tau)$    contains a non-trivial  converging sequence then
{\it asdim} $ ( G, \mathcal{S}_{\tau}) = \infty$.

\vskip 6pt

Proof.} We denote
$\mathcal{F}=\{\mathcal{S}_{\tau(a_{n})}: (a_{n}) $   is a  non-trivial
 sequence  in $(G,\tau )$  converging to  $0$ in  $(G, \tau )\}$.
Let  $A\in \mathcal{S}_{\tau}$.
By the  definition of $\mathcal{S}_{\tau}$,   there exists a subset
 $B\in  \mathcal{S}_{\tau}$,
  $A \subseteq B$   and a finite number
  $(a_{1n})_{n\in\omega}    , \ldots  , (a_{mn})_{n\in\omega}$
   of sequences converging to  $0$  in
$(G, \tau)$  such that $B$  can be obtained
 from the sets
 $\{ a_{1n}, - a_{1n}: n\in\omega\}, \ldots  ,
  \{ a_{mn}, - a_{mn}: n\in\omega\}$
 by the finite number  of additions  of these sets and join of finite subsets of $G$.
We choose a sequence
$(b_{n})_{n\in \omega}$ converging  to $0$  in $(G, \tau)$  and containing each
$(a_{1n})_{n\in\omega}    , \ldots  , (a_{mn})_{n\in\omega}$
 as a subsequence.
Then $A\in   \mathcal{S}_{\tau(b_{n})} $,
$\mathcal{S}_{\tau(b_{n})} \subseteq \mathcal{S}_{\tau}$.
Apply Theorems 3 and 4. \hfill   $\Box$
 \vskip 7pt

We note that Theorem 5 answers Question 2 from \cite{b3}.
 \vskip 6pt

{\bf Question 1. } {\it  Let  $(G,\tau )$ be a countable non-discrete metrizable abelian group.
Does $(G, \mathcal{S}_{\tau} )$ contain an asymorphic copy  of $\mathbb{H}$?}
\vskip 8pt

In \cite{b4},  the authors ask about asymptotic dimension of
$(G, \mathcal{C}_{\tau} )$,  where $G$  is an infinite cyclic subgroup of the circle.
We put this question in more genera form.
\vskip 6pt

 {\bf Question 2. } {\it  Let  $(G,\tau )$ be a countable non-discrete metrizable group.
 Is  {\it asdim} $ ( G, \mathcal{C}_{\tau}) = \infty$? }
\vskip 7pt

 {\bf Example 2. } Let  $G$ be the direct sum of $\omega$ copies  of
 $\mathbb{Z}_{2}$
   endowed with the topology induced by the  Tikhonov topology of
   $\mathbb{Z}_{2}^{\omega}$.
Since
{\it asdim} $ ( G, \mathcal{C}_{\tau}) > 0$
and
$(G, \mathcal{C}_{\tau})$
 is asymorphic to
 $(G, \mathcal{C}_{\tau})\times (G, \mathcal{C}_{\tau})$,
  we  see that
  {\it asdim} $ ( G, \mathcal{C}_{\tau}) = \infty$.
\vskip 5pt

Let $(X,  \mathcal{E})$  be a coarse space. A function $f: X \longrightarrow  \{0,1\}$  is called {\it slowly oscillating} if,  for every $E\in  \mathcal{E}$ ,  there exists a bounded subset $B$ of $X$  such that $f|_{E[x]}  =  const $  for each  $x\in X\setminus   B$.
We endow $X$  with the discrete  topology,  identify the Stone-$\check{C}$ech  compactification $\beta X$  of $X$  with the set of ultrafilters on $X$   and denote
 $X^{\sharp} = \{ p\in  \beta X : $  each $P\in p$  is unbounded $\}$.
We define an equivalence $\sim$  on  $X^{\sharp}$  by the rule:  $p\sim q$
 if and only if   $f^{\beta} (p) = f^{\beta} (q) $ for every slowly oscillating
 function $f: X\longrightarrow\{0,1\} $.
The quotient $X^{\sharp} / \sim$  is called a space of  {\it ends }  or {\it binary corona}  of   $(X, \mathcal{E})$,   see  [12, Chapter 8].

\vspace{5 mm}

{\bf Theorem 6}. {\it If a non-trivial $T$-sequences    $(a_{n}) _{n \in\omega} $
generates $G$ then the space of  ends of
$( G, \mathcal{C}_{\tau(a_{n}) })$ is a singleton.

\vspace{5 mm}

Proof.} First we show that for every slowly oscillating function $f: G\longrightarrow\{0,1\}$ there exists an $m$ such that
$$
f|_{G\setminus  A _{m}}  =  const .\eqno(5).
$$

Indeed, by the definition of slow oscillation and Theorem 1,  there exists $m\in \omega$  such that $f|_{x+A}  =  const $ for each $x\in G\setminus A_{m}$.
We show nw that (5) holds true for this $m$.

We take arbitrary  $y, z \in  G\setminus  A _{m} $. Since $A$ generates $G$ and contains 0, there exists
an index $k$ such that $y,z \in A_k$, i.e.,
$$
y= b_{1} + \ldots  + b_{k} \ \ \mbox{ and }\  \ z= c_{1} + \ldots + c_{k},
$$
for appropriate
   $ \  \  b_{1}, \ldots ,b_{k}$,   $ \ \  c_{1}, \ldots ,  c_{k} \in  A $.
 By a property of $T$-sequences established at the end of the proof of Theorem 3,   there exists a member $a_{m_1}$ of $(a_n)$  such that
 $$
 a_{m_1} +  y  \notin  A _{m+1} \ \ \mbox{ and } \ \ \  a_{m_1} +  z  \notin  A _{m+1},
 $$
 since $(a_n)$ is a $T$-sequence.
Then
$a_{m_1} +    b_{2} + \ldots  +  b_{k} \notin  A _{m}$,
   $ \ \  a_{m_1} +  c_{2}+  \ldots  +   c_{k} \notin  A _{m}$.
   Therefore,
$$
 f( a_{m_1} +  b_{2} + \ldots  +  b_{k}) =  f (y) \ \mbox{ and }\  \  f( a_{m_1} +  c_{2}+  \ldots  +   c_{k} )=  f(z). \eqno(6)
$$
Repeating this trick $k$   times, we can replace
$b_{2}, \ldots  ,  b_k$ and $c_{2} , \ldots  ,  c_k$, by appropriate  members $a_{m_2}, \ldots , a_{m_k}$ of $(a_n)$, as before.
 Hence, we can replace $b_{2} + \ldots  +  b_k$ and $c_{2} + \ldots  +  c_k$, by $a_{m_2} + \ldots  +  a_{m_k}$ in (6).
This obviously gives $f(y) = f(z)$ and proves (5).

 Finally, to prove the assertion of the theorem, pick $p,q\in  X^{\sharp}$. In order to check that $p\sim q$ fix an arbitrary
slowly oscillating function $f: X\longrightarrow\{0,1\} $. We have prove that $f^{\beta} (p) = f^{\beta} (q) $. Pick an $m$ with (5). Since $G\setminus  A _{m}\in p \cap q$, for every $P\in p$ and for every $Q\in q$ we have     $P_1:= P \setminus A _{m}\in p $ and $Q_1:= Q\setminus  A _{m}\in q$
 and $f|_{P_1} = f|_{Q_1}$ is constant in view of (5). This proves that     $f^{\beta} (p) = f^{\beta} (q) $.
\hfill  $ \Box$
\vspace{6 mm}

{\bf Question 3. } {\it  Let  $(G,\tau )$ be a countable non-discrete metrizable abelian group.
Is the space of ends of
$(G, \mathcal{S}_{\tau} )$  a singleton?
The same question for}
$(G, \mathcal{C}_{\tau} )$.
\vskip 10pt

\section{Metrizability and normality}

By  [12, Theorem 2.1.1],  a coarse space $(X, \mathcal{E})$  is metrizable if and only if $\mathcal{E}$  has a countable base.

\vspace{6 mm}

{\bf Theorem 7}. {\it  For an infinite  abelian group $G$  and a $T$-sequence  $(a_{n}) _{n \in\omega}$  in $G$,
 $(G, \mathcal{S} _{\tau(a_{n})})$ is metrizable  if and only if $G$  is countable.
\vspace{6 mm}

Proof.}  Apply Theorem 1.  \hfill   $\Box$

\vspace{6 mm}

{\bf Theorem 8}. {\it   Let $G, $  be a non-discrete  metrizable group.
Then the coarse structure $\mathcal{E}_{\tau}$  of  $(G, S)$ does not have a linearly  ordered base.
In particular,  is  not metrizable.

\vspace{6 mm}
Proof.}
We assume that $\mathcal{E}_{\tau}$  has a linear base $\mathcal{E}$  and choose a sequence
 $(A_{n})_{n\in\omega} $  in $\mathcal{E}$  such that
 $A_{n} \subset A_{n+1}$ and the closure of
 $\bigcup_{n\in\omega}A_{n}$
   is not  compact.
If
 $B\in \mathcal{E}$   then
 $B\subseteq  A_{m}$  for some  $m$  so
  $\{ A_{n}:  n\in\omega\}$  is a base for $\mathcal{E}_{\tau}$.

Now let $\{ U _{n}: n\in \omega\}$  is a base  of neighbourhoods  of the  identity of  $(G, \tau)$.
We choose an injective sequence $(b_{n})_{n\in \omega}$  in $G$   such that   $b_{n}\in U_{n}$,   $b_{n}\notin A_{n}$.
Then
$\{ b_{n} : n\in \omega \} \in \mathcal{S}_{\tau}$
 but
 $\{ b_{n} : n\in \omega \} \setminus A_{n} \neq\emptyset$.
Hence,
$\{ A_{n} : n\in \omega \}$
is not a base for $\mathcal{E}_{\tau}$  and we get a contradiction. \hfill   $\Box$

\vspace{7 mm}

Let $(X, \mathcal{E})$ be a ballean.
A subset $U$ of $X$ is called an {\it asymptotic  neighbourhood} of a subset
 $Y\subseteq  X$  if,  for  every  $E\in \mathcal{E}$,  $E[Y]\setminus  U$ is bounded.
\vspace{7 mm}

Two subsets  $Y, Z$  of $X$  are called
\vspace{7 mm}

\begin{itemize}
\item{}  {\it  asymptotically disjoint  if, for every $E\in \mathcal{E}$,
 $E[Y]\cap   E[Z]$  is bounded;
 \vskip 7pt

\item{}     asymptotically separated if $Y, Z$ have disjoint asymptotic neighbourhoods.}
\end{itemize}
\vskip 7pt

A ballean $(X, \mathcal{E})$  is called {\it normal}  [8] if any two
 asymptotically disjoint subsets of  $X$  are asymptotically separated.
  Every  ballean with linearly ordered  base is normal  [8, Proposition 1.1].

We suppose that a non-discrete metrizable group  $(G, \tau)$ is  topologically isomorphic  to the product
 $G_{1} \times G_{2}$  of  infinite groups.
Applying Theorem 8 and  Theorem 1.4 from [1],  we  conclude that  $(G, \mathcal{S}_{\tau})$  is not  normal.

\vskip 6pt

{\bf Question 4. } {\it  Let $(G, \tau)$  be a non-discrete metrizable group.
Is the ballean  $(G, \mathcal{S}_{\tau})$    non-normal?}

\vskip 10pt

\section{Functional boundedness}

Let  $(X, \mathcal{E})$ be a coarse space.
Following [2],  we say that a function
$f: (X, \mathcal{E}) \longrightarrow \mathbb{R}$  is

 \vspace{7 mm}

\begin{itemize}
\item{}  {\it   bornologous} if  $f(B)$ is bounded in $\mathbb{R}$ for each bounded  subset  $B$  of $X$;
\vspace{4 mm}

\item{}  {\it    macro-uniform} if,  for every
$E\in\mathcal{E}$,
  the  supremum
   $\sup _{x\in X} \ diam \ f(E[x])$
   is finite;

\vspace{4 mm}

\item{}  {\it   eventually  macro-uniform} if, for every
$E\in\mathcal{E}$,
  there exists a bounded  subset  $B $  of  $X$  such that
$\sup _{x\in X\setminus B}  \ diam  f(E[x])$
   is finite;

\vspace{4 mm}

\item{}  {\it  slowly oscillating} if, for  every  $E\in\mathcal{E}$    and   $\varepsilon>0$, there  exists a
bounded  subset  $B$ of $X$  such  that

$ diam  f(E[x])< \varepsilon$   for each  $x\in X\setminus B$.

\end{itemize}
\vskip 7pt

We say that a  coarse space $(X, \mathcal{E})$
 is

  \vspace{7 mm}

\begin{itemize}
\item{}  {\it  $b$-bounded } if each  bornologous  function on $X$  is bounded;
\vspace{4 mm}

\item{} {\it mu-bounded}  if each macro-uniform   function on $X$  is bounded;
\vspace{4 mm}

 \item{} {\it  emu- bounded} if,  for  every   macro-uniform  function  $f$ on $X$,   there exists a bounded   subset   $B$  of $X$  such that  $f$  is bounded on  $X\setminus  B$;
  \vspace{4 mm}

\item{} {\it   so-bounded} if, for every  slowly   oscillating  function $f$ on $X$,  there exists a bounded subset $B$  of $X$  such that $f$  is bounded on  $X\setminus  B$.

\end{itemize}
\vskip 7pt

If $(G, \tau)$ is discrete then  every  function
$f: (G,  \mathcal{S}_{\tau})\longrightarrow\mathbb{R}$
 is  bornologous.  If
 $(G, \tau)$
   has a non-trivial  converging sequence then  there is a  non-bornologous  function
   $f: (G,  \mathcal{S}_{\tau})\longrightarrow\mathbb{R}$.

\vspace{6 mm}

{\bf Theorem 9}. {\it  If
$(G, \tau)$
 is  metrizable  and compact then
 $(G, \mathcal{S}_{\tau})$
   is  $b$-bounded.
If $(G, \tau)$
 is countable  and  metrizable   then
 $(G, \mathcal{S}_{\tau})$
   is not $b$-bounded.

\vspace{6 mm}
Proof.} In the first case, we  assume  that  there  exists an  unbounded bornologous function
$f: (G,  \mathcal{S}_{\tau})\longrightarrow\mathbb{R}$.

We choose a sequence
$(a_{n})_{n\in\omega}$
 in  $G$  such  that
 $|f(a_{n})|> n$.
Passing to a subsequence,  we may  suppose that
$(a_{n})_{n\in\omega}$
  converges  to  some  point $a$.
Then $f$ is not bounded on the set
$\{ a, a_{n}: n\in\omega   \}  \in\mathcal{S}_{\tau}$
 so  $f$  is not  bornologous.

In the second case, we   denote by $H$   the completion of $G$  and  choose
 $a n$  injective  sequence
 $(b_{n})_{n\in\omega}$
  in  $G$   converging to some point
   $h\in H\setminus  G$.
Then  we define a  function
$f: (G, \tau) \longrightarrow\mathbb{R}$
 by $f(b_{n})= n$ and $f(x)=0$ for each
 $x\in G\setminus \{ b_{n} : n\in\omega\}$.
Since
$K\cap\{ b_{n} : n\in\omega \}$
is finite for each compact subset $K$  of $G$,  we see that $f$  is bornologous.
Hence,  $(G, \mathcal{S}_{\tau})$ is not  $b$-bounded.\hfill   $\Box$

\vspace{7 mm}

If $(G, \tau)$ is  discrete and
$(G, \mathcal{S}_{\tau})$
 is  $m u$-bounded then $G$  is a Bergman group [2].

\vspace{6 mm}

{\bf Theorem 10}. {\it  If
$(G, \tau)$
 is  metrizable  and
totally bounded then
$(G, \mathcal{S}_{\tau})$
 is $mu$-bounded.

\vspace{6 mm}
Proof.}
We assume the contrary.  Then there exists a macro-uniform  function
$f: (G,  \mathcal{S}_{\tau})\longrightarrow\mathbb{R}^{+}$
  such that
  $f(a_{n+1}) - f(a_{n})> n$
   for some sequence
   $(a_{n})_{n\in\omega}$
   in $(G, \tau)$.
Passing to a subsequence, we may suppose that
$(a_{n})_{n\in\omega}$
 converges to some point $h\in H$,  where $H$  is the completion of $G$.
We  denote $b_{n} = a_{n} a _{n+1}  ^{-1}$    and  observe that the sequence
$(b_{n})_{n\in\omega}$
 converges to the identity $e$  of $G$.
We put
$B=\{ e, b_{n}: n\in\omega \}$ and note that
 $a_{n}\in B  a _{n+1}$.
Since $B\in \mathcal{S}_{\tau}$,  we see that $f$  is  not macro-uniform.

\vskip 6pt

{\bf Question 5. } {\it Let be a compact metrizable abelian group. Is $(G, \mathcal{S}_{\tau})$ $ \ so$-bounded?
$ \ mu$-bounded?}

\vspace{6 mm}

CONTACT INFORMATION

I.~Protasov: \\
Faculty of Computer Science and Cybernetics  \\
        Kyiv University  \\
         Academic Glushkov pr. 4d  \\
         03680 Kyiv, Ukraine \\ i.v.protasov@gmail.com

\end{document}